\documentclass[a4paper,12pt]{article}
    \usepackage[top=2.5cm,bottom=2.5cm,left=2.5cm,right=2.5cm]{geometry}
    \usepackage{cite, amsmath, amssymb}
\usepackage{tikz}
\usepackage{amsmath}

\begin{document}
\begin{center}
{\LARGE\bf Laplacian eigenvalues of equivalent cographs}
\end{center}
\begin{center}
{\large \bf J. Lazzarin, O.F. M\'arquez and F. C. Tura}
\end{center}
\begin{center}
\it Departamento de Matem\'atica, UFSM, Santa Maria, RS, 97105-900, Brazil
\end{center}
\begin{center}
{\small lazzarin@smail.ufsm.br, oscar.f.marquez-sosa@ufsm.br, fernando.tura@ufsm.br}
\end{center}


\newcommand{\bfx}{{\mathbf x}}
\newcommand{\casei}{{\bf case~1}}
\newcommand{\subia}{{\bf subcase~1a}}
\newcommand{\subib}{{\bf subcase~1b}}
\newcommand{\subic}{{\bf subcase~1c}}
\newcommand{\caseii}{{\bf case~2}}
\newcommand{\subiia}{{\bf subcase~2a}}
\newcommand{\subiib}{{\bf subcase~2b}}
\newcommand{\caseiii}{{\bf case~3}}
\newcommand{\myvar}{x}

\newtheorem{Thr}{Theorem}
\newtheorem{Pro}{Proposition}
\newtheorem{Que}{Question}
\newtheorem{Con}{Conjecture}
\newtheorem{Cor}{Corollary}
\newtheorem{Lem}{Lemma}
\newtheorem{Fac}{Fact}
\newtheorem{Ex}{Example}
\newtheorem{Def}{Definition}
\newtheorem{Prop}{Proposition}
\def\floor#1{\left\lfloor{#1}\right\rfloor}

\newenvironment{my_enumerate}{
\begin{enumerate}
  \setlength{\baselineskip}{14pt}
  \setlength{\parskip}{0pt}
  \setlength{\parsep}{0pt}}{\end{enumerate}
}

\newenvironment{my_description}{
\begin{description}
  \setlength{\baselineskip}{14pt}
  \setlength{\parskip}{0pt}
  \setlength{\parsep}{0pt}}{\end{description}
}


\begin{abstract} Let $G$ and $H$ be  equivalent cographs with their reduction  $R_G$ and $R_H,$  and suppose the vertices of $R_G$ and $R_H$ are labeled by the twin numbers $t_i $ of the $k$ twin classes they represent. In this paper,  we prove that $G$ and $H$ have at least $k + \sum_{i\in I}(t_i-1)$ Laplacian eigenvalues in common, where  $I$ is the indices of the twin classes whose types are identical in $G$ and $H.$ This confirms the conjecture proposed by T. Abrishami \cite{Abris}.
We also show that no two nonisomorphic equivalent cographs are $L$-cospectral.
\end{abstract}
$\hspace{1cm}${\it keywords:}  Laplacian eigenvalues, twins numbers,  cograph, $L$-cospectral graphs

$\hspace{0.25cm}$ {\it AMS subject classification:} 15A18, 05C50, 05C85.
\baselineskip=0.30in


\section{Introduction}
\label{intro}

$\hspace{0,5cm}$ A cograph is a simple graph which contains no path on four vertices as induced subgraph, namely it is $P_4$-free graph.
An equivalent definition (see \cite{Stewart}) is that cographs can be obtained recursively by the following rules: $(i)$ a graph on a single vertex is a cograph,
  $(ii)$ a finite  union and join of two cographs are  cographs.
 This allow us to represent this class of graphs through an unique rooted tree $T_G,$  called the cotree. For more details, see Section \ref{Sec2}.

An important subclass of cographs are the threshold graphs. 
Threshold graphs also can be defined in terms of forbidden subgraphs, namely they are $\{P_4, 2K_2, C_4\}$-free graphs. For an account on different characterizations and properties of threshold graphs, one can see 
\cite{Mah95} and the references therein.

Let $G=(V,E)$ be a graph, we say two vertices $v$ and $w$ are {\em twins} if $N(v)-w=N(w)-v,$ where $N(v)$ denotes the neighborhood of $v.$ A twin partition of a graph $G$ is the partition of the vertices into their equivalence classes under  the relation  of being  twins,  denoted by $ V(G)= T_1 \cup T_2 \cup \ldots \cup T_k.$ The {\em twins numbers} $t_1, t_2, \ldots, t_k$ of a graph $G$ are the size of twins classes.

\begin{figure}[h!]
       \begin{minipage}[c]{0.25 \linewidth}
\begin{tikzpicture}[ultra thick]
  [scale=1,auto=left,every node/.style={circle,scale=0.7}]
  \node[draw,circle,fill=gray!30!white] (p) at (-3,5) {};
  \node[draw,circle,fill=green] (o) at (0,6) {};
  \node[draw,circle,fill=green] (n) at (-2,6) {};
   \node[draw,circle,fill=gray!30!white] (q) at (-2,3) {};
    
  \node[draw,circle,fill=red] (f) at  (1,5) {};
    \node[draw,circle,fill=red] (f1) at  (1,4) {};

  \node[draw,circle,fill=white] (c) at (0,3) {};
  \node[draw,circle,fill=gray!30!white] (d) at (-3,4) {};
  \path

        (n) edge node[below]{}(p)
        (o) edge node[below]{}(p)
        (q) edge node[below]{}(p)
       (d) edge node[below]{}(p)
       (d) edge node[below]{}(q)
       (o) edge node[below]{}(q)

        (o) edge node[below]{}(d)
        (n) edge node[below]{}(d)
        
        (f) edge node[below]{}(c)
        (f1) edge node[below]{}(c)
        
        (f) edge node[below]{}(o)
        (f1) edge node[below]{}(o)
        
        (f1) edge node[below]{}(q)
        (f) edge node[below]{}(q)

        (f) edge node[below]{}(p)
        (f1) edge node[below]{}(p)

        (f) edge node[below]{}(d)
        (f1) edge node[below]{}(d)
        
        (f) edge node[below]{}(f1)
        
        (f1) edge node[below]{}(n)
        
        (n) edge node[below]{}(f)
        (n) edge node[below]{}(q);
\end{tikzpicture}
       \end{minipage}\hfill
       \begin{minipage}[l]{0.4 \linewidth}
\begin{tikzpicture}[ultra thick]
\draw[fill=green ](0,0) circle [radius=0.7];
\draw[fill](0, 0.25) circle[radius =0.1];
\draw[fill](0, -0.25) circle[radius =0.1];

\draw[fill=white ](0,3) circle [radius=0.7];
\draw[fill](0, 3) circle[radius =0.1];

\draw(0.7, 3)--(2.3,3);
\draw(3, 2.3)--(3,0.7);
\draw(0.7,0)--(2.3,0);

\draw(0.5, 0.5)--(2.5,2.5);


\draw[fill=red](3,3) circle[radius=0.7];
\draw(3, 3.25)--(3,2.75);
\draw[fill](3, 3.25) circle[radius =0.1];
\draw[fill](3, 2.75) circle[radius =0.1];

\draw[fill=gray!30!white ](3,0) circle [radius=0.7];
\draw[fill](3, 0.25) circle[radius =0.1];
\draw[fill](2.75, -0.25) circle[radius =0.1];
\draw[fill](3.25, -0.25) circle[radius =0.1];
\draw(3, 0.25)--(2.75,-0.25);
\draw(3, 0.25)--(3.25,-0.25);
\draw(2.75, -0.25)--(3.25,-0.25);

\end{tikzpicture}
       \end{minipage}
       \caption{A cograph and its reduction representation}
       \label{fig1}
\end{figure}
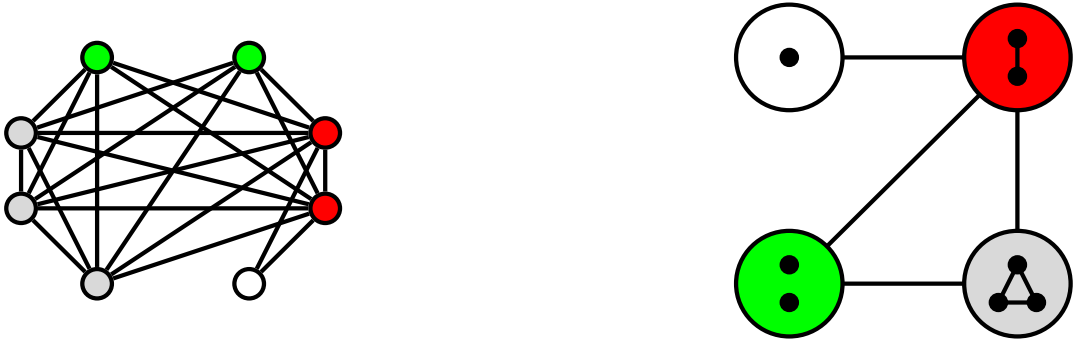

Let $G$ be a cograph with twin classes   $T_1, T_2, \ldots, T_k.$ The {\em twin reduction}  of graph $G,$ denoted  $R_{G},$ is the subgraph induced  by $\{u_1, \ldots, u_k\},$ where $u_i \in T_i$  is a representative  of class $T_i.$ The Figure \ref{fig1} shows a cograph $G$ and its twin reduction  $R_G.$ The twins classes $T_1, T_2, T_3,$ and $T_4,$ where $T_1$ is the green vertices, $T_2$ is the gray vertices, $T_3$ is the white vertex and $T_4$ is the red vertices. We say two cographs $G$ and $H$ are equivalent, if their reduction $R_G$ and $R_H$ are isomorphic and if the twin numbers of the vertices $V(R_G)$ and $V(R_H)$ are identical.

Our motivation for considering cographs comes from spectral graph theory. There is a considerable body of knowledge on the spectral properties of cographs and threshold graphs
related to adjacency matrix \cite{Allem, Allem2,Milica,Cesar2,Bapat,BSS2011,gorbani,JTT2016,JOT2019,lou,Moha,Royle,Tura21}.
However, the literature does not seem to provide many articles about the Laplacian matrices of cographs. One of those sporadic works and very well known is the paper of Russel Merris \cite{Merris} which shows that the nonzero Laplacian eigenvalues of threshold graphs are equal to Ferrer's conjugate of its degree sequence.

A recent and interesting work about the Laplacian eigenvalues of cographs  is the Master's thesis presented by Tara Abrishami \cite{Abris}. In this work, a characterization of cograph Laplacian eigenvalues is given. In particular, it was proposed the following conjecture:
\begin{Con}\label{conjec}
Let $G$ and $H$ be  equivalent cographs with their reduction   $R_G$ and $R_H,$ and suppose the vertices of $R_G$ and $R_H$ are labeled by the twin numbers $t_i$ of the $k$ twin classes they represent.  Then $G$ and $H$ have at least $k + \sum_{i\in I}(t_i-1)$ Laplacian eigenvalues in common, where  $I$ is the indices of the twin classes whose types are identical in $G$ and $H.$
\end{Con}

In this paper, we prove that conjecture holds. As an  immediate result, we  show that no two equivalent cographs $G$ and $H$ are $L$-cospectral graphs. The main tool used to prove  the conjecture,  and some  known results are reviewed  in  Section \ref{Sec2}.  In Section \ref{Sec3}, we  characterize equivalent cographs in terms of their cotrees. In Section \ref{Sec4}, we confirm that conjecture holds. In the final section, we show how to find the $L$-cospectral linear size families of cographs, from a pair of two nonisomorphic $L$-cospectral cographs.

\section{Notations and Preliminaries} \label{Sec2}

$\hspace{0,5cm}$ Let  $G= (V,E)$ be an undirected graph with vertex set $V$ and edge set $E,$ without loops or multiple edges.  For $v\in V,$  $N(v)$ denotes the  {\em open neighborhood}  of $v,$ that is, $\{w|\{v,w\}\in E\}.$  The  {\em closed neighborhood}  $N[v] = N(v) \cup \{v\}.$ If $|V| = n,$ the   Laplacian matrix $L(G)$ of a graph $G$ is given by $L(G)=\delta(G) -  A(G),$ where $\delta(G)$ is the degree matrix of $G$ and $A(G)$ is the adjacency matrix of $G.$ A value $\mu(G)$ is a Laplacian eigenvalue of $G$
 if $\det(L(G) - \mu I_n ) = 0$, and since
$L(G)$ is real symmetric and  positive defined, the Laplacian eigenvalues of $G$ are real numbers non-negative.

\subsection{Cotrees}

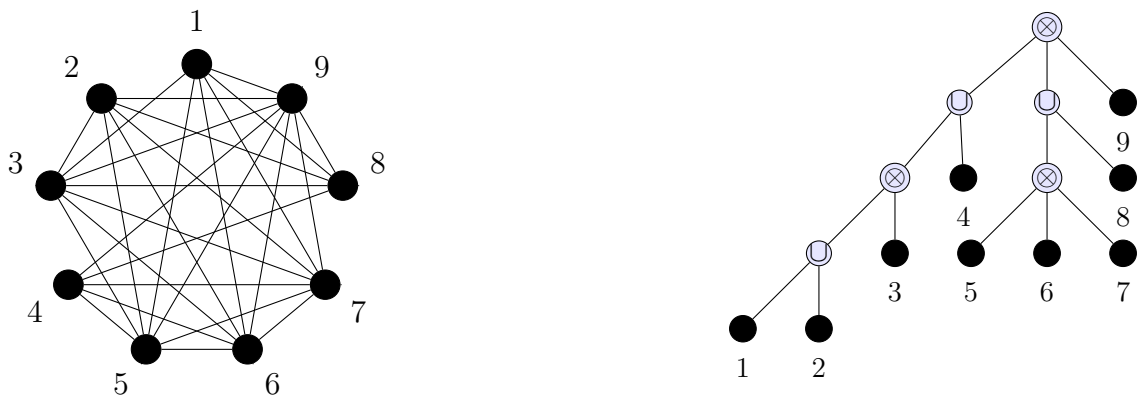
\begin{figure}[h!]
       \begin{minipage}[c]{0.25 \linewidth}
\begin{tikzpicture}
  [scale=0.65,auto=left,every node/.style={circle}]
  \foreach \i/\w in {1/,2/,3/,4/,5/,6/,7/,8/,9/}{
    \node[draw,circle,fill=black,label={360/9 * (\i - 1)+90}:\i] (\i) at ({360/9 * (\i - 1)+90}:3) {\w};} 
  \foreach \from in {3,5,6,7,9}{
    \foreach \to in {1,2,...,\from}
      \draw (\from) -- (\to);}
       \foreach \from in {8}{
       \foreach \to in {1,2,3,4,\from}
       \draw (\from) -- (\to)
       ;}
\end{tikzpicture}
       \end{minipage}\hfill
       \begin{minipage}[l]{0.4 \linewidth}
\begin{tikzpicture}
  [scale=1,auto=left,every node/.style={circle,scale=0.9}]
  \node[draw,circle,fill=black,label=below:$1$] (p) at (-2,3) {};
  \node[draw,circle,fill=black,label=below:$7$] (o) at (3,4) {};
  \node[draw,circle,fill=black,label=below:$6$] (n) at (2,4) {};
   \node[draw,circle,fill=black,label=below:$5$] (q) at (1,4) {};
  \node[draw, circle, fill=blue!10, inner sep=0] (m) at (2,5) {$\otimes$};
  \node[draw,circle,fill=black,label=below:$8$] (l) at (3,5) {};
  \node[draw,circle,fill=black,label=below:$9$] (k) at (3,6) {};
  \node[draw, circle, fill=blue!10, inner sep=0] (j) at (2,6) {$\cup$};
  \node[draw,circle,fill=blue!10, inner sep=0] (h) at (2,7) {$\otimes$};
  \node[draw, circle, fill=blue!10, inner sep=0] (g) at (0.85,6) {$\cup$};
  \node[draw,circle,fill=black,label=below:$4$] (f) at  (0.9,5) {};
  \node[draw,circle,fill=blue!10, inner sep=0] (a) at (0,5) {$\otimes$};
  \node[draw, circle, fill=blue!10, inner sep=0] (b) at (-1,4) {$\cup$};
  \node[draw,circle,fill=black,label=below:$3$] (c) at (0,4) {};
  \node[draw,circle,fill=black,label=below:$2$] (d) at (-1,3) {};
  \path (a) edge node[left]{} (b)
        (a) edge node[below]{} (c)
        (b) edge node[left]{} (d)
        (f) edge node[right]{}(g)
        (g) edge node[left]{}(a)
        (h) edge node[right]{}(j)
        (h) edge node[left]{}(g)
        (h) edge node[left]{}(k)
        (j) edge node[right]{}(l)
        (j) edge node[below]{}(m)
        (m) edge node[below]{}(n)
        (m) edge node[right]{}(o)
        (b) edge node[left]{} (p)
        (m) edge node[left]{} (q);
\end{tikzpicture}
       \end{minipage}
       \caption{A cograph $G=((((v_{1}\cup v_{2})\otimes v_{3})\cup v_{4})\otimes(((v_{5}\otimes v_{6})\otimes v_{7})\cup v_{8}))\otimes v_{9}$. and its cotree.}
       \label{fig2}
\end{figure}

A cotree $T_G$ of a cograph $G$ is a rooted tree in which any interior vertex $w$ is either of $\cup$ type (corresponds to union) or $\otimes$ type (corresponds to join). The terminal vertices (leaves) are typeless and represent the vertices of the cograph $G.$  We say that {\em depth} of the cotree is the number of edges of  the longest path from the root to a leaf. To build a cotree for a connected cograph, we simply place a $\otimes$ at the tree's root, placing $\cup$ on interior vertices with odd depth, and placing $\otimes$ on interior vertices with even depth. 
 All interior vertices have at least two children. In \cite{BSS2011} this structure is called {\em minimal  cotree}, but  throughout this paper we call it simply a cotree.
 The  Figure \ref{fig2} shows  a cograph and its cotree with depth equals  to 4.

Two vertices $u$ and $v$  are {\em duplicate}  if $N(u) = N(v)$ and {\em coduplicate} if $N[u]=N[v].$  In fact, any collection of mutually coduplicate (resp. duplicate) vertices, e.g. with the same neighbors and adjacent (resp. not adjacent),  have a common parent of type $ \otimes$ (resp.  $\cup$).

\noindent{\bf Remark:} We note that  vertices in a twin class $T_i$ of $G$ correspond to coduplicate (resp. duplicate) vertices in $T_G,$ if they are pairwise adjacent (resp. nonadjacent) twins.

\subsection{Diagonalization}\label{}

For comparing the Laplacian eigenvalues of two cographs $G$ and $H,$ we use a straightforward translation of an algorithm due to Jacobs et al. \cite{JTT2016}  to the context of Laplacian matrices  of cographs. The original  algorithm    constructs a {\em diagonal} matrix congruent to $A + x I_n$,
where $A$ is the adjacency matrix of a cograph,
and $x$ is an arbitrary scalar, using $O(n)$ time and space.

One of the advantages of this method is that it can be slightly  modified  in such a way  that we can determine, for any $-x \in \mathbb{R},$ the number of Laplacian eigenvalues of a cograph $G$ that are larger  than $x,$ equal to $x$ and smaller than $x,$ respectively.  The algorithm's input  is the cotree $T_G$ and $x \in \mathbb{R}.$
Each leaf $v_i, i =1, \ldots,n$ have a value $d_i$ that represents the diagonal element of  $L(G)+xI_n.$ It initializes  all entries $d_i$ with $\delta(v_i) +x,$ where $\delta(v_i)$ denotes the degree of vertex $v_i.$
Each iteration, a pair $\{v_k, v_l\} $ of the  duplicate or coduplicate vertices with maximum depth  is selected.
 Then they are processed, that is,  assignments are given to $d_k$ and $d_l,$ such that either one or both rows (columns) are diagonalized. When a $k$ row(column)  corresponding to vertex  $v_k$  has been diagonalized then $v_k$ is  removed from the $T_G,$ it means that $d_k$ has a permanent final value. Then the algorithm moves to the cotree $T_G -v_k.$  The algorithm is shown in Figure~\ref{algo}.

It is worth to mention that for each iteration,  the algorithm executes one of the six subcases.
It should be noted that \subia~ and \subiia~ are the normal cases,
and the other four subcases represent singularities.
Executing \subib~ requires $\beta = -1$,
executing \subiib~ requires $\beta = 0$,
executing \subic~ requires $\alpha + \beta = -2,$
and executing {\bf subcase 2c} requires $\alpha + \beta = 0$.

\begin{figure}[h]
{\tt
\begin{tabbing}
aaa\=aaa\=aaa\=aaa\=aaa\=aaa\=aaa\=aaa\= \kill
     \> INPUT:  cotree $T_G$, scalar $\myvar$\\
     \> OUTPUT: diagonal matrix $D=[d_1, d_2, \ldots, d_n]$ congruent to $L(G) + \myvar I_n$\\
     \>\\
     \>   $\mbox{ Algorithm}$ Diagonal $(T_{G}, x)$ \\
     \> \> initialize $d_i := \delta(v_i)+ \myvar$, for $ 1 \leq i \leq n$ \\
     \> \> {\bf while } $T_G$  has $\geq 2$    leaves      \\
     \> \> \>  select a pair $(v_k, v_l)$  (co)duplicate of maximum depth with  parent $w$\\
     \> \> \>     $\alpha \leftarrow  d_k$    $\beta \leftarrow d_{l}$\\
     \> \> \> {\bf if} $ w=\otimes$\\
     \> \> \> \>  {\bf if} $\alpha + \beta \neq -2$  \verb+                //subcase 1a+    \\
     \> \> \> \> \>   $d_{l} \leftarrow \frac{\alpha \beta -1}{\alpha + \beta +2};$ \hspace*{0,25cm} $d_{k} \leftarrow \alpha + \beta +2; $\hspace{0,25cm}   $T_G = T_G - v_k$ \\
     \> \> \> \>  {\bf else if } $\beta=-1$ \verb+                //subcase 1b+   \\
     \> \> \> \> \>   $d_{l} \leftarrow -1$ \hspace*{0,25cm}   $d_k  \leftarrow 0;$ \hspace{0,25cm} $T_G = T_G - v_k$ \\
     \> \> \> \>  {\bf else  }  \verb+                      //subcase 1c+   \\
     \> \> \> \> \>   $d_{l} \leftarrow -1$  \hspace*{0,25cm} $d_k \leftarrow (1+\beta)^2;$ \hspace{0,25cm} $T_G= T_G -v_k;$ \hspace{0,25cm} $T_G = T_G -v_l$  \\
     \> \> \>     {\bf else if} $w=\cup$\\
     \> \> \> \>  {\bf if} $\alpha + \beta \neq 0$  \verb+               //subcase 2a+    \\
     \> \> \> \> \>   $d_{l} \leftarrow \frac{\alpha \beta}{\alpha +\beta};$ \hspace*{0,25cm}   $d_k \leftarrow \alpha +\beta;$ \hspace{0,25cm} $T_G = T_G - v_k$ \\
     \> \> \> \>  {\bf else if } $\beta=0$ \verb+                //subcase 2b+   \\
     \> \> \> \> \>   $d_{l} \leftarrow 0;$ \hspace*{0,25cm}  $d_k  \leftarrow 0;$ \hspace{0,25cm} $T_G = T_G - v_k$ \\
     \> \> \> \>  {\bf else  }  \verb+                      //subcase 2c+   \\
     \> \> \> \> \>   $d_{l} \leftarrow \beta;$  \hspace*{0,25cm} $v_k \leftarrow -\beta;$ \hspace{0,25cm} $T_G =T_G - v_k;$ \hspace{0,25cm} $T_G = T_G - v_l$  \\

     \> \>  {\bf end loop}\\
\end{tabbing}
}
\caption{\label{algo} Diagonalization algorithm}
\end{figure}

Now, we will present a few results which the proofs are similar to work  \cite{JTT2016}. The following theorem is based on Sylvester's Law of Inertia. 

\begin{Thr}
\label{main1}
Let $G$ be a cograph and let $(d_v)_{v \in T_G}$ be the sequence produced by Diagonalize $(T_G,-x).$
Then the diagonal matrix $D = diag(d_v)_{v \in T_G}$ is congruent to $L(G) +xI_n,$ so that  the number of (positive - negative - zero) entries in $(d_v)_{v\in T_G}$ 
is equal to the number eigenvalues of $L(G)$  that are (greater than $\myvar$ - small than $\myvar$ - equal to $\myvar$).  
\end{Thr}

The following two lemmas show that if a vertex $\otimes$ or $\cup$, in the cotree, have leaves with the same value, then, we can use the following routines.

\begin{Lem}
\label{lem1}
If $v_1, \ldots, v_m$ have parent $w= \otimes$, each with the same diagonal value $y \neq -1$, then the algorithm performs $m-1$ iterations of  \subia~ assigning, during iteration  $j:$
\begin{equation}
d_k  \leftarrow \frac{j+1}{j}(y+1) \hspace{0,5cm} d_l  \leftarrow \frac{y-(j-1)}{j+1}
\end{equation}
\end{Lem}

\begin{Lem}
\label{lem2}
If $v_1, \ldots, v_m$ have parent $w= \cup$, each with the same diagonal value $y\neq 0$, then the algorithm performs $m-1$ iterations of  \subiia~ assigning, during iteration  $j:$
\begin{equation}
d_k  \leftarrow \frac{j+1}{j}y \hspace{0,5cm} d_l  \leftarrow \frac{y}{j+1}
\end{equation}
\end{Lem}

\section{Equivalent Cographs and their cotrees}\label{Sec3}

$\hspace{0,5cm}$ In this section, we characterize equivalent cographs in terms of their cotrees. Let $G$ and $H$ be equivalent cographs with their reduction $R_G$ and $R_H.$ If $R_G \cong R_H$ and fixed $T_G,$ we will show how to get $T_H$ from $T_G.$

Let $G$ be a cograph and $T_G$ its cotree. Let $u, v$ be leaves in the cotree  $T_G$ which have the lowest common ancestor an interior vertex, represented by $\mbox{lca}(u,v).$ Clearly, they are adjacent if and only if  $\mbox{lca}(u,v)= \otimes$. The following result can be verified immediatly.

\begin{Lem}\label{lem3}
Let $G$ and $H$ be equivalent cographs with their cotrees $T_G$ and $T_H.$
Let $u,v$ be leaves in $T_G$ which are neither coduplicate nor duplicate vertices, and let $u', v'$ be their corresponding leaves in $T_H.$ 
Then $\mbox{lca}(u,v)$ and $\mbox{lca}(u',v')$ are the same type.
\end{Lem}

\begin{Def}
Let $G$ be a cograph and $T_G$ its cotree. For any pair $u,v \in T_G,$ we define the distance between $u$ and $v$ in $T_G,$ denoted by $dist_{T_G}(u,v),$ as the shortest path of  interior vertices between them. 
\end{Def}

\begin{Def}
Let $G$ and $H$ be equivalent cographs with their reduction $R_G$ and $R_H.$
Let $u$ be a representative of the twin class $T_u \in G,$ and let $u'$ be  its corresponding in $T_{u'} \in H.$ We say $u=u',$ if $T_u$ and $T_{u'}$ are twin classes of same type. Otherwise, we say  $u\neq u'.$
\end{Def}

\begin{Lem}
\label{lem4}
Let $G$ and $H$ be equivalent cographs with their cotrees $T_G$ and $T_H.$
Let $u, v$ be the  representatives of the twin classes $T_u, T_v \in G,$ and let $u', v'$ be  their respective correspondents in $T_{u'}, T_{v'} \in H.$
\begin{enumerate}
 \item[(i)] If $u=u'$ and $v=v'$ then  $dist_{T_G}(u,v) = dist_{T_H}(u',v').$
\item[(ii)]   If $u=u'$ and $v\neq v'$ then  $dist_{T_G}(u,v) = dist_{T_H}(u',v') \pm 1.$
\item[(iii)]   If $u\neq u'$ and $v\neq v'$ then  $dist_{T_G}(u,v) = dist_{T_H}(u',v') \pm 2.$

\end{enumerate}
\end{Lem}
\noindent {\bf Proof:}
We prove the item $(i).$ We assume that $lca(u,v)=lca(u',v') =\otimes.$ By contradiction, we suppose that $dist_{T_G}(u,v) < dist_{T_H}(u',v').$ Since that $u=u'$ and $v=v',$ we have that
\begin{equation}\label{eq3}
dist_{T_H}(u',v') = dist_{T_G}(u,v) + 2l
\end{equation}
for some positive integer $l \geq 1.$

\begin{figure}[h!]
\begin{center}
\begin{tikzpicture}
  [scale=0.85,auto=left,every node/.style={circle,scale=0.85}]
  
  
  \node[draw,circle,fill=black] (n) at (4,4) {};
  \node[draw,circle,fill=black] (n1) at (11,4) {};
  
   \node[draw,circle,fill=blue!10,inner sep=0] (q) at (3,4) {$\otimes$};
    \node[draw,circle,fill=blue!10,inner sep=0] (qt1) at (10,4) {$\otimes$};
    \node[draw,circle,fill=black, label=below:$u'$] (r1) at (10,3) {};

    \node[draw,circle,fill=black, label=below:$u$] (q1) at (3,3) {};

  \node[draw, circle, fill=blue!10, inner sep=0] (m) at (4,5) {$\cup$};
   \node[draw, circle, fill=blue!10, inner sep=0] (m1) at (11,5) {$\cup$};
 
  
  \node[draw,circle,fill=blue!10, inner sep=0] (k2) at (6,5) {$\cup$};
\node[draw,circle,fill=black, label=below:$v$] (k3) at (6,4) {};
   \node[draw,circle,fill=black, label=below:$t'$] (k1) at (13,4) {};
 \node[draw,circle,fill=blue!10, inner sep=0] (k4) at (7,4) {$\otimes$};
  \node[draw,circle,fill=black] (k5) at (7,3) {};

  \node[draw, circle, fill=blue!10, inner sep=0, label=right:$w_i$] (j) at (5,6) {$\otimes$};
   \node[draw, circle, fill=blue!10, inner sep=0, label=right:$w'_i$] (j1) at (12,6) {$\otimes$};
 \node[draw,circle,fill=black] (jk) at (12,5) {};

 \node[draw,circle,fill=blue!10, inner sep=0] (h1) at (13,5) {$\cup$};
 \node[draw,circle,fill=blue!10, inner sep=0] (h2) at (14,4) {$\otimes$};
\node[draw,circle,fill=black, label=below:$w'$] (l1) at (14,3) {};
   \node[draw,circle,fill=blue!10, inner sep=0] (h3) at (15,3) {$\cup$};
   \node[draw,circle,fill=black, label=below:$v'$] (r2) at (15,2) {};

   \node[draw,circle,fill=black] (a) at (5,5) {}; 
  
  
  
   \path 
        (jk) edge node[right]{}(j1)
        (h1) edge node[right]{}(k1)
 (h1) edge node[right]{}(h2)

        (h1) edge [dashed] node[right]{}(j1)
 (k2) edge node[right]{}(k3)
         (k2) edge node[right]{}(k4)
      (k5) edge node[right]{}(k4)

        (h2) edge [dashed]  node[left]{}(h3)
        (j1) edge [dashed] node[right]{}(m1)
        
        (j) edge[dashed] node[right]{}(k2)
         (j) edge node[right]{}(a)
        (j) edge [dashed] node[below]{}(m)
        (l1) edge node[below]{}(h2)
        (q) edge node[below]{}(q1)
        (m) edge node[below]{}(n)
         (m1) edge node[right]{}(n1) 
          (m1) edge node[right]{}(qt1)
           (r1) edge node[right]{}(qt1)
 (r2) edge node[right]{}(h3)
        (m) edge node[left]{} (q);

\end{tikzpicture}
       \caption{The partial cotrees $T_G$ and $T_H.$}
       \label{fig4}
\end{center}
\end{figure}
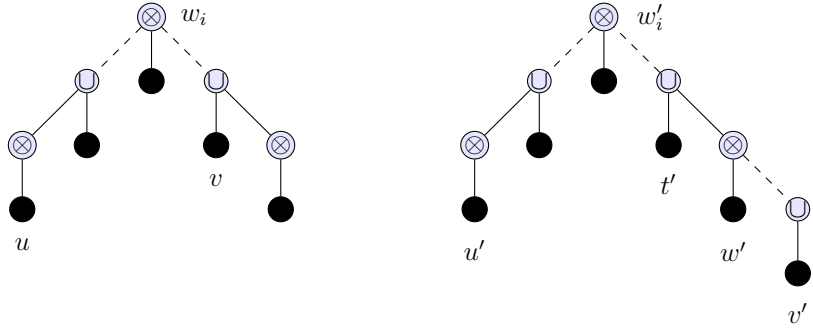

Without loss of generality, we assume that $u$ and $u'$ are in  twin classes of type $\otimes,$ while that $v$ and $v'$ are in twin classes of type $\cup,$ and their partial cotrees $T_G$ and $T_H$ are represented in the Figure  \ref{fig4}.

Now consider the respective reduction graphs $R_G$ and $R_H$ of cographs $G$ and $H.$ From equation (\ref{eq3}) follows that there are vertices $t', w' \in R_H$ such that $w'\sim v'$ and $t'\nsim v'.$
Since $R_G$ and $R_H$ are isomorphic graphs, there are vertices $t,w \in R_G$  with the same properties.

We claim the leaf $t, w\in T_G$ are in the same branch that leaf $v.$ 
Since $t\sim u$ and for any leaf  in a different branch that $v$ implies being adjacent to $v,$ follows the statement. Now, if $w$ is in a different branch that $v$ it implies that $w \sim t,$ what contradics $w\nsim t.$ 
If the leaf $t$ is below to $v$ and $w\sim v,$ it implies $w\sim t,$ a contradicition. Now, if the leaf $t$ is above to $v,$ and since $w\sim v, w\nsim t,$ it implies that $w$ is between $t$ and $v,$ and therefore we must have $l=0,$ in the equation (\ref{eq3}). If $lca(u,v)=lca(u',v') =\cup,$ the proof is analogous.

The proof is similar for the items $(ii)$ and $(iii).$ $\hspace{6.5cm}\square$

Given a cograph $G,$ we note that its reduction $R_G$ is obtained by taking only one representative of each twin classes of $G.$
In terms of cotree, it means, if we remove all exceed leaves of $T_G,$ we have a cotree which represents $T_{R_G}.$ 
If $R_G$ and $R_H$  are isomorphic graphs by Lemma \ref{lem4} we have $|dist_{T_{G}}(u,v) - dist_{T_{H}}(u',v')| \leq 2,$ for any leaves $u,v \in T_G$ and their corresponding
$u',v' \in T_H.$ This allows us to claim that $T_H$ can be obtained from $T_{G},$  as the following result:

\begin{Thr}
\label{main2}
Let $G$ and $H$ be cographs with their reduction $R_G$ and $R_H.$ If $T_{G}$ and $T_{H}$ are the cotrees of $G$ and $H$, then $R_G \cong R_H$  if and only if  fixed $T_{G}$ and for some interior vertex $w_i  \in T_{G},$ having leaves $t_i \geq 2,$   then $T_{H}$ is obtained from $T_{G}$ by one of following operations: 
\begin{enumerate}
 \item[(i)] adding an interior vertex,  one level below to $w_i,$ and taking their leaves $t_i.$
\item[(ii)]  removing the  interior vertex $w_i$ which has no interior vertex as successor, whose father has no leaves and taking their leaves $t_i.$ 
\end{enumerate}
\end{Thr}

\noindent {\bf Proof:} Let $G$ and $H$ be  equivalent cographs with their reductions $R_G$ and $R_H.$ If $T_{G}$ and $T_{H}$ are the cotrees of  $G$ and $H,$ respectively,   according to Lemma \ref{lem4}, for each pair of leaves $u$ and $v$ in $T_{G}$ which are neither coduplicate nor duplicate vertices and their corresponding $u'$ and $v'$ in $T_{H},$ we have
\begin{equation}
 |dist_{T_{G}}(u,v) - dist_{T_{H}}(u',v')| \leq 2 
\end{equation}

If the distance is preserved  and since the $lca(u,v)$ and $lca(u',v')$ are the same type in both cotrees then  $T_{G}$ and $T_{H}$ are the same. If the distance increased or decreased by one, and taking into account that the $lca(u,v) =w_i$ and $lca(u',v') = w'_i$ are the same type, we have two possibilities: was creating a new interior vertex below to $w_i$ and taking their leaves $t_i$ of $w_i,$ or was removing the vertex $w_i,$ which has no interior vertex as successor, whose father has no leaves and taking their leaves $t_i.$  Finally, if the distance increase or decrease by two, then either of operations $(i)$ or $(ii)$ occurs twice.

Now, let $T_G$ and $T_H$ be the cotrees of equivalent cographs $G$ and $H,$ respectively. We just need to check that reductions $R_G$ and $R_H$ are isomorphic.
Let $T_G$ and $T_H$ be the cotrees of $G$ and $H.$ We assume that $T_H$ is obtained from $T_G$ under the operations $(i)$ and $(ii).$ First, we note that the number of twin classes are preserved, since that the only operation allowed is to become a coduplicate vertices into duplicate vertices or vice versa. Second, from the operations $(i)$ and $(ii)$
the $lca(u,v)$ in $T_G$ and its corresponding $lca(u',v')$ in $T_H$ are the same type, which implies that the adjacencies of $R_G$ and $R_H$ are preserved. 
Therefore, thus $R_G$ and $R_H$ are isomorphic. $\hspace{7,5cm}\square$

\section{The proof of conjecture}\label{Sec4}
$\hspace{0,5cm}$ The next lemmas will be used to prove the main results of this section:
\begin{Lem}
\label{lem5} Let $G$ and $H$ be equivalent cographs. Let $u$ be a representative of twin class $T_u$ of $G$ with twin number $t_u,$ and let $u'$ be its corresponding in a twin class $T_{u'}$ of $H.$
If $\delta(u)$ denotes the degree of vertex $u$ then
$$\delta(u') = \left\{
\begin{array}{lr}
\delta(u) & \mbox{if $u'=u$   }\\
 \delta(u) +(t_u -1) & \mbox{if $u\neq u'$ and $T_u$ is a clique set.}
\end{array} \right.$$
\end{Lem}

\noindent{\bf Proof:}
Let $u$ and $u'$ be the  representatives of twin classes $T_u$ of $G$ and $T_{u'}$ of $H,$ respectively. If $u=u'$ then obvious we have that $\delta(u')=\delta(u).$
Now, we assume that $u\neq u',$ and $T_{u}$ is a clique set. Taking into account that  vertex $u'$  will be disconnected only of the $t_{u} -1$ vertices of same class $T_{u'},$ follows that $\delta(u')=\delta(u) +(t_u-1),$ as desired.$\hspace{1,5cm}\square$

\begin{Lem}
\label{lem6} Let $G$ be a cograph with twin classes $T_1, T_2, \ldots, T_k,$ and twin numbers $t_1, t_2, \ldots, t_k.$ Let $u_i$ be a representative of twin class $T_i.$ If $\delta(u_i)$ denotes the degree of $u_i,$ then
$$\mu(G) = \left\{
\begin{array}{lr}
\delta(u_i) & \mbox{if $T_{u_i}$ is a coclique set}\\
 \delta(u_i) + 1 & \mbox{if $T_{u_i}$ is a clique set}
\end{array} \right.$$
 is a Laplacian eigenvalue of $G$ with multiplicity at least $t_i -1,$ for $i=1,2,\ldots,k.$
\end{Lem}

\noindent{\bf Proof:}
Let $G$ be a cograph and let $u_i$ be a representative of twin class $T_{u_i}$ with twin number $t_i,$ for $i=1,2,\ldots,k.$ Now, we considere Diagonalization of $(T_G,x),$ with
$x=-\delta(u_i),$ if $T_{u_i}$ is a coclique set and $x=-\delta(u_i)-1,$ if $T_{u_i}$ is a clique set, for $i=1,2,\ldots, k.$  Since that coduplicate (respect. duplicate) vertices of $T_G$ correspond to  clique (respect. coclique) set in $G,$ then after initialization $(d_i = \delta(u_i)+x)$, we have the following values for the leaves of $T_G$
$$\left\{
\begin{array}{lr}
-1 & \mbox{for coduplicate vertices}\\
 0 & \mbox{ for duplicate vertices.}
\end{array} \right.$$

From this, it is easy to see that for coduplicate vertices the \subib~occurs, while that for the duplicate vertices the \subiib~occurs. In both cases, the algorithm assigns a zero as a permanent value. Since each twin class $T_i$ have $t_i$ vertices, follows for each iteration we have at least $t_i-1$ zeros, as desired. $\hspace{7,5cm}\square$

\begin{Thr}
Let $G$ and $H$ be  equivalent cographs with their reduction   $R_G$ and $R_H,$ and suppose the vertices of $R_G$ and $R_H$ are labeled by the twin numbers $t_i$ of the $k$ twin classes they represent.  Then $G$ and $H$ have at least $k + \sum_{i\in I}(t_i-1)$ Laplacian eigenvalues in common, where  $I$ is the indices of the twin classes whose types are identical in $G$ and $H.$
\end{Thr}

\noindent{\bf Proof:}
Let $G$ and $H$ be  equivalent cographs with their reduction   $R_G$ and $R_H,$ and suppose the vertices of $R_G$ and $R_H$ are labeled by the twin numbers $t_i$ of the $k$ twin classes they represent.  Let $I$ be the indices of the twin classes whose types are identical in $G$ and $H$ 
with cardinality $ 0 \leq  | I | < k.$

In order for proving the conjecture, we note that a cograph $G$ of order $n=\sum t_i,$ 
 each Laplacian eigenvalue $\mu(G)$ of $G$ belongs one of following subsets
\begin{equation}
\sum_{i \in I} (t_i -1) \cup (\bigcup_{i=1}^{ | I |} t_i) \cup (\bigcup_{i=1}^{k-| I |} t_i)
\end{equation}

We first will show that $G$ and $H$ have $\sum_{i\in I}(t_i -1)$ Laplacian eigenvalues in common.
Let $T_G$ and $T_H$ be the cotrees of $G$ and $H,$ respectively. They have $\sum_{i \in I} (t_i -1) \cup (\bigcup_{i=1}^{ | I |} t_i)  $ leaves in common. Since $R_G \cong R_H,$ and $u_i, u'_i $ are the respective  representatives of twin classes which are identical in $G$ and $H$  having the same degree $\delta(u_i)= \delta(u'_i),$ hence by Lemma \ref{lem6}, we have at least
\begin{equation}
\sum_{i \in I} (t_i -1)
\end{equation}
Laplacian eigenvalues in common.

Now, let $\mu(G)$ be a Laplacian eigenvalue of $G$ which $\mu(G) \in (\bigcup_{i=1}^{ | I |} t_i) \cup  (\bigcup_{i=1}^{k-| I |} t_i) .$  We will show that $\mu(G)$ is also one of the $k$ Laplacian eigenvalues of $H.$  For this, we apply the Diagonalization algorithm simultaneously in both cotrees $T_G$ and $T_H$ with $x=-\mu(G).$
It is sufficient to show when the algorithm assigns a zero in Diag$(T_G,x),$ we also must have a zero in Diag$(T_H,x).$

Let $\alpha_i$ and $\alpha'_i$ be the assigments given in the $i$-th iteration during execution of Diag$(T_G,x)$ and Diag$(T_H,x),$ respectively. Obviously, we have $\alpha_i = \alpha'_i,$ if $T_G$ and $T_H$ are identical. It remains to be seen when $T_G$ and $T_H$ have different types of cotrees but $R_G$ and $R_H$ are isomorphic.

Suppose that $T_G$ has an interior vertex $w_i$ having $t_i \geq 2$ coduplicate vertices and a pendant vertex with assigment $\alpha_i,$
while that $T_H$ has an interior vertex $w'_i$ having no leaves but the same pendant vertex with same assigment $\alpha_i$ and an interior vertex as successor having $t_i \geq 2$ duplicate vertices, as the Figure \ref{fig6} has shown.

\begin{figure}[h!]
\begin{center}
\begin{tikzpicture}
  [scale=0.85,auto=left,every node/.style={circle,scale=0.85}]
  
  
  \node[draw,circle,fill=black, label=below:$\alpha_i$] (n) at (1,4) {};
  \node[draw,circle,fill=black,  label=below:$\alpha_i$] (n1) at (8,4) {};
  


  \node[draw, circle, fill=blue!10, inner sep=0] (m) at (1,5) {$\cup$};
   \node[draw, circle, fill=blue!10, inner sep=0] (m1) at (8,5) {$\cup$};
 
  \node[draw,circle,fill=black,label=below:$1$] (l1) at (9,4) {};
  
  \node[draw,circle,fill=black, ,label=below:$t_i$] (k) at (3,5) {};
   \node[draw,circle,fill=black,label=below:$t_i$] (k1) at (10,4) {};
   
  \node[draw, circle, fill=blue!10, inner sep=0,label=right:$w_i$] (j) at (2,6) {$\otimes$};
   \node[draw, circle, fill=blue!10, inner sep=0, label=right:$w'_i$] (j1) at (9,6) {$\otimes$};
%

 \node[draw,circle,fill=blue!10, inner sep=0] (h1) at (10,5) {$\cup$};
   
   \node[draw,circle,fill=black, ,label=below:$1$] (a) at (2,5) {}; 
  
  
  
   \path 
        (h1) edge node[right]{}(k1)
        
        (h1) edge node[right]{}(j1)
        
        (j1) edge node[right]{}(m1)
        
        (j) edge node[right]{}(k)
         (j) edge node[right]{}(a)
        (j) edge node[below]{}(m)
        (l1) edge node[below]{}(h1)
        (m) edge node[below]{}(n)
       (m1) edge node[right]{}(n1); 

\end{tikzpicture}
       \caption{The partial cotrees $T_G$ and $T_H.$}
       \label{fig6}
\end{center}
\end{figure}
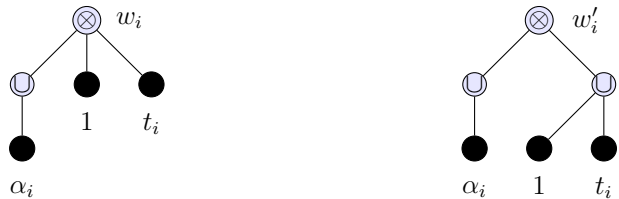

Let $v_i$ and  $v'_i $ be the respective  representatives of twin classes which are not identical in $G$ and $H.$ Obviously, by Lemma \ref{lem6}, we can assume that $\mu(G)$ differs of  $\delta(v_i)+1$ and $\delta(v'_i).$
Applying the algorithm in the coduplicate vertices of $T_G,$ since all $t_i$  leaves have the same value $y = \delta(v_i) +x,$ by Lemma \ref{lem1}, after $t_i -1$ iterations, we have a pendant vertex with value 
\begin{equation}\label{eq6}
d_l = \frac{y -(t_i -1)}{t_i +1} = \frac{\delta(v_i) -\mu(G) -(t_i -1)}{t_i +1}
\end{equation}

Now, we apply the algorithm in the duplicate vertices of $T_H,$ since all $t_i$  leaves have the same value $y' = \delta(v'_i) +x,$ by Lemma \ref{lem2}, after $t_i -1$ iterations, we have a pendant vertex with value 
\begin{equation}\label{eq7}
d'_l = \frac{y'}{t_i +1} = \frac{\delta(v'_i) -\mu(G)}{t_i +1}
\end{equation}

We claim that $d_l = d'_l.$ From equations (\ref{eq6}) and (\ref{eq7}), follows 
$$  \frac{\delta(v_i) -\mu(G) -(t_i -1)}{t_i +1} =  \frac{\delta(v'_i) -\mu(G)}{t_i +1} \Leftrightarrow \delta(v_i) -\mu(G) -(t_i -1)= \delta(v'_i) -\mu(G)$$
\begin{equation}\label{eq8}
\delta(v'_i) = \delta(v_i)  + (t_i -1)
\end{equation}
which accords with the Lemma \ref{lem5}. This shows that the algorithm will assign the same value to both $T_G$ and $T_H.$
Since $-x=\mu(G)$ is a Laplacian eigenvalue of $G$ and a zero should be assigned during execution of Diag$(T_G,x)$ then
either situations can occurs: a zero was assigned in the previous step to $d_l,$ which correponds a leaf that belongs to $ (\bigcup_{i=1}^{ | I |} t_i),$ or a zero is assigned exactly after to process the values $\alpha_i$ and $d_l.$
Thus, a zero must be assigned during execution of Diag$(T_H,x).$

The proof is similar if the $lca(v_i, v'_i)$ is of  $\cup$ type. 
Therefore, thus follows the result as desired. $\hspace{14cm}\square$

\begin{Cor}
Let $G$ and $H$ be equivalent cographs with their reduction $R_G$ and $R_H$ having $k\geq 2$ vertices.  If $G$ and $H$ are $L$-cospectral graphs
then $G \cong H.$
\end{Cor}

\noindent{\bf Proof:}
Let $G$ and $H$ equivalent cographs with their reduction $R_G$ and $R_H.$ We proceed by induction on the number $k\geq 2$ of vertices of $R_G$ and $R_H.$ 
The base case, $k=2$ is trivial to verify.

We assume that the result holds for any two equivalent cographs $G$ and $H$ with their reduction $R_G$ and $R_H$ having $k-1$ vertices.
Now let $G'$ and $H'$ be equivalent cographs $L$-cospectral with their reduction $R_{G'}$ and $R_{H'}$ having $k$ vertices.

We note that there is a vertex $u\in R_{G'}$ and its respective corresponding $v\in R_{H'}$ such that $\delta(u)=\delta(v),$ and $T_u$ and $T_v$ are twin classes of same type.  Otherwise, we  have two distinct Laplacian eigenvalues, according to Lemma \ref{lem6}.
Then, the cographs $G'-T_u \cong G$ and $H'-T_v \cong H$ are $L$-cospectral graphs. By the induction hypothesis we have $G\cong H$ and therefore follows that $G' \cong H'.$ $\hspace{9cm}\square$

\section{$L$-cospectral cographs}\label{Sec5}

$\hspace{0,5cm}$ Two nonisomorphic graphs with the same $L$-spectrum are called $L$-cospectral. In this last section, we show how to find the $L$-cospectral linear size families of cographs, from a pair of two nonisomorphic $L$-cospectral cographs.

The next lemma is very well known and it will be used  for our construction:

\begin{Lem}
\label{lem7} Let $G$ and $H$  be a graphs on  $n_1$ and $n_2$ vertices, respectively. If $0=\mu_1(G) \leq \mu_2(G) \leq \ldots \leq \mu_1(G)$ and  $0=\mu_1(H) \leq \mu_2(H) \leq \ldots \leq \mu_1(H)$
are the  Laplacian eigenvalue of $G$  and $H,$ respectively. Then the Laplacian eigenvalues of $ G \otimes H$ are
$$ 0, n_2 + \mu_2(G), n_2 + \mu_3(G), \ldots, n_2 + \mu_1(G),$$
$$  n_1 + \mu_2(H), n_1 + \mu_3(H), \ldots, n_1 + \mu_1(H), n_1 + n_2.$$
\end{Lem}

For each integer $n\geq 3,$ we define the following cographs of order $2n +1$
\begin{itemize}
\item  $G_{2n+1} = n K_1 \otimes (K_n \cup K_1);$
\item  $H_{2n+1}= (((n -1) K_1 \otimes K_1)\cup K_1) \otimes (K_{n-1} \cup K_1).$
\end{itemize}

The Figure \ref{fig6} shows the reduction representation of $G_{9}$ and $H_9.$

\begin{figure}[h!]
       \begin{minipage}[c]{0.25 \linewidth}
\begin{tikzpicture}[ultra thick]
\draw[fill=white ](0,0) circle [radius=0.7];
\draw[fill](0, 0) circle[radius =0.1];

\draw[fill=white ](0,3) circle [radius=0.7];
\draw[fill](0.25, 3.25) circle[radius =0.1];
\draw[fill](0.25, 2.65) circle[radius =0.1];
\draw[fill](-0.25, 3.25) circle[radius =0.1];
\draw[fill](-0.25, 2.65) circle[radius =0.1];

\draw(0, 0.7)--(0,2.3);
\draw(0.7, 3)--(2.3,3);



\draw[fill=gray!30!white](3,3) circle[radius=0.7];
\draw(2.75, 3.25)--(2.75,2.75);
\draw(3.25, 3.25)--(3.25,2.75);
\draw(2.75, 3.25)--(3.25,3.25);
\draw(3.25, 2.75)--(2.75,2.75);
\draw(3.25, 3.25)--(2.75,2.75);
\draw(2.75, 3.25)--(3.25,2.75);
\draw[fill](2.75, 3.25) circle[radius =0.1];
\draw[fill](2.75, 2.75) circle[radius =0.1];
\draw[fill](3.25, 3.25) circle[radius =0.1];
\draw[fill](3.25, 2.75) circle[radius =0.1];

\end{tikzpicture}
       \end{minipage}\hfill
       \begin{minipage}[l]{0.4 \linewidth}
\begin{tikzpicture}[ultra thick]
\draw[fill=white ](0,0) circle [radius=0.7];
\draw[fill](0, 0) circle[radius =0.1];

\draw[fill=white ](0,-2.5) circle [radius=0.7];
\draw[fill](0, -2.5) circle[radius =0.1];
\draw(0.7, -2.2)--(2.6,2.4);

\draw(0.7, -2.2)--(2.4,0);

\draw[fill=white ](0,3) circle [radius=0.7];
\draw[fill](0.25, 2.75) circle[radius =0.1];
\draw[fill](-0.25, 2.75) circle[radius =0.1];
\draw[fill](0, 3.25) circle[radius =0.1];

\draw(0, 0.7)--(0,2.3);
\draw(0.7, 3)--(2.3,3);
\draw(0.7,0)--(2.3,0);
\draw(0.7, 3)--(3,0.7);
\draw(0.5, 0.5)--(2.5,2.5);


\draw[fill=gray!30!white](3,3) circle[radius=0.7];
\draw(3, 3.25)--(2.75,2.75);
\draw(3, 3.25)--(3.25,2.75);
\draw(2.75, 2.75)--(3.25,2.75);
\draw[fill](3, 3.25) circle[radius =0.1];
\draw[fill](2.75, 2.75) circle[radius =0.1];
\draw[fill](3.25, 2.75) circle[radius =0.1];

\draw[fill=gray!30!white ](3,0) circle [radius=0.7];
\draw[fill](3, 0) circle[radius =0.1];

\end{tikzpicture}
       \end{minipage}
       \caption{The reduction representation of $G_{9}$ and $H_{9}$}
       \label{fig6}
\end{figure}
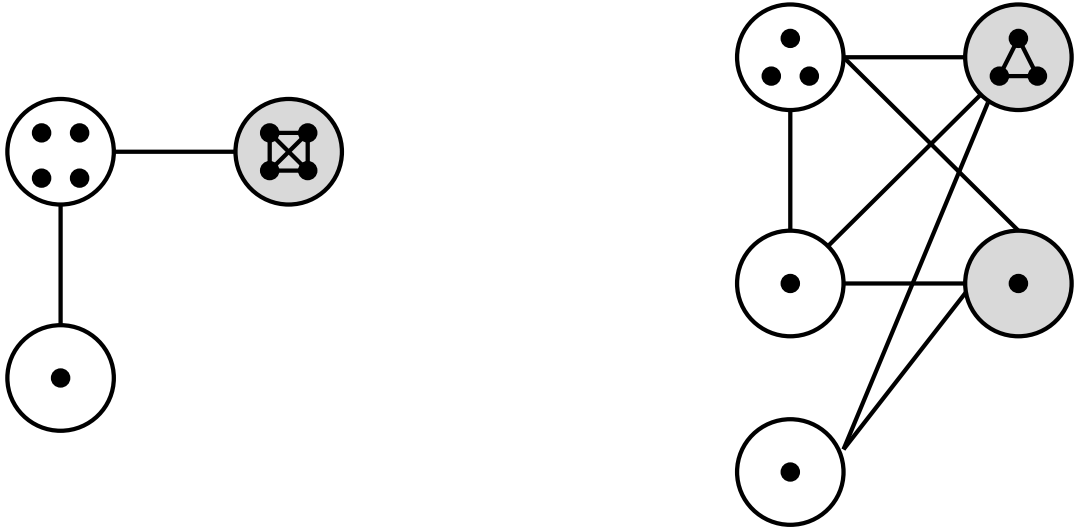

\begin{Lem}
\label{lem8}
The cographs $G_{2n+1}$ and $H_{2n+1}$ of order $2n+1$ defined above are nonisomorphic and $L$-cospectral graphs.
\end{Lem}

\noindent{\bf Proof:} It is obvious that $G_{2n+1}$ and $H_{2n+1}$ are nonisomorphic graphs. Since the $L$-eigenvalues of $K_n \cup K_1$ are $n$ and $0$ with multiplicity $n-1$ and $2,$
by Lemma \ref{lem7}, we have that the $L$-spectrum of $G_{2n+1}$ is 
$$ 0 + (n+1), 0 +( n+1) , \ldots, 0+(n+1),  0+n, 0 + n, \ldots, n+n, (n+1) +n.$$
Therefore, the $L$-spectrum of $G_{2n+1}$ is $2n+1, 2n, n+1,n,0$ with their respective multiplicities $1, n-1, n-1, 1, 1.$  By similar calculus, we have $H_{2n+1}$ have the same $L$-spectrum.$\square$


\begin{Thr}
Let $G^\prime$ be a cograph of order $n.$ Then $G^\prime \otimes G_{2n+1}$ and $ G^\prime \otimes H_{2n+1} $  are nonisomorphic and $L$-cospectral graphs.
\end{Thr}

\noindent{\bf Proof:} Let $G^\prime$ be a cograph of order $n.$ Let $G_{2n+1}$ be a cograph of order $2n+1.$ 
If $ 0=\mu_1 \leq \mu_2 \leq \ldots \leq \mu_1$ are the $L$-eigenvalues of $G^\prime,$ by Lemma \ref{lem7} the $L$-eigenvalues of $G^\prime \otimes G_{2n+1}$ are
$$ \mu_1 +(2n+1), \mu_2 +(2n+1), \ldots, \mu_1 +(2n+1),$$
$$  n+ n, n+(n+1), \ldots, n+(2n), n+(2n+1), n+(2n+1).$$

Since $G_{2n+1}$ and $H_{2n+1}$ are $L$-cospectral graphs, then the $L$-eigenvalues of $ G^\prime \otimes H_{2n+1} $ are obtained by same procedure above.
Therefore, thus $G^\prime \otimes G_{2n+1}$ and $ G^\prime \otimes H_{2n+1} $  are nonisomorphic and $L$-cospectral graphs. $\hspace{9cm} \square$

\end{document}